\magnification=1200
\input amssym.def
\input amssym.tex

\hsize=13.5truecm
\baselineskip=16truept
\baselineskip=16truept
\font\secbf=cmb10 scaled 1200
\font\eightrm=cmr8
\font\sixrm=cmr6

\font\eighti=cmmi8

\font\sixi=cmmi6
\skewchar\eighti='177 \skewchar\sixi='177

\font\eightsy=cmsy8
\font\sixsy=cmsy6
\skewchar\eightsy='60 \skewchar\sixsy='60

\font\eightit=cmti8

\font\eightbf=cmbx8
\font\sixbf=cmbx6

\let\sc=\tensc

\font\eightsc=cmcsc10 scaled 800
\font\secbf=cmb10 scaled 1200
\font\subsecfont=cmb10 scaled \magstephalf
\font\amb=cmmib10

\font\ambi=cmmib10 scaled 700

\newfam\mbfam 

\textfont\mbfam\amb \scriptfont\mbfam\ambi


\def\aa{\def\rm{\fam0\eightrm}%
  \textfont0=\eightrm \scriptfont0=\sixrm \scriptscriptfont0=\fiverm
  \textfont1=\eighti \scriptfont1=\sixi \scriptscriptfont1=\fivei
  \textfont2=\eightsy \scriptfont2=\sixsy \scriptscriptfont2=\fivesy
  \textfont3=\tenex \scriptfont3=\tenex \scriptscriptfont3=\tenex
  \def\sc{\eightsc}
  \def\it{\fam\itfam\eightit}%
  \textfont\itfam=\eightit
  \def\bf{\fam\bffam\eightbf}%
  \textfont\bffam=\eightbf \scriptfont\bffam=\sixbf
   \scriptscriptfont\bffam=\fivebf
  \normalbaselineskip=9.7pt
  \setbox\strutbox=\hbox{\vrule height7pt depth2.6pt width0pt}%
  \normalbaselines\rm}

\def\Proof{\vskip12pt\noindent{\bf Proof.} }

\def\Def#1{\vskip12pt\noindent{\bf Definition #1}}
\def\Remark#1{\vskip12pt\noindent{\bf Remark #1}}

\def\m@th{\mathsurround=0pt}

\def\cc#1{\hbox to .89\hsize{$\displaystyle\hfil{#1}\hfil$}\cr}
\def\lc#1{\hbox to .89\hsize{$\displaystyle{#1}\hfill$}\cr}
\def\rc#1{\hbox to .89\hsize{$\displaystyle\hfill{#1}$}\cr}

\def\eqal#1{\null\,\vcenter{\openup\jot\m@th
  \ialign{\strut\hfil$\displaystyle{##}$&&$\displaystyle{{}##}$\hfil
      \crcr#1\crcr}}\,}

\def\section#1{\vskip 22pt plus6pt minus2pt\penalty-400
        {{\secbf
        \noindent#1\rightskip=0pt plus 1fill\par}}
        \par\vskip 12pt plus5pt minus 2pt
        \penalty 1000}

\def\subsection#1{\vskip 20pt plus6pt minus2pt\penalty-400
        {{\subsecfont
        \noindent#1\rightskip=0pt plus 1fill\par}}
        \par\vskip 8pt plus5pt minus 2pt
        \penalty 1000}

\def\subsubsection#1{\vskip 18pt plus6pt minus2pt\penalty-400
        {{\subsecfont
        \noindent#1}}
        \par\vskip 7pt plus5pt minus 2pt
        \penalty 1000}

\def\center#1{{\begingroup \leftskip=0pt plus 1fil\rightskip=\leftskip
\parfillskip=0pt \spaceskip=.3333em \xspaceskip=.5em \pretolerance 9999
\tolerance 9999 \parindent 0pt \hyphenpenalty 9999 \exhyphenpenalty 9999
\par #1\par\endgroup}}

\def\\{\hfill\break}

\def\kwadrat{\null\ \hfill\null\ \hfill$\square$}
\def\mida#1{{{\null\kern-4.2pt\left\bracevert\vbox to 6pt{}\!\hbox{$#1$}\!\right\bracevert\!\!}}}
\def\midy#1{{{\null\kern-4.2pt\left\bracevert\!\!\hbox{$\scriptstyle{#1}$}\!\!\right\bracevert\!\!}}}

\def\diagint{{\raise1.5pt\hbox{$\scriptscriptstyle\diagup$}\hskip-8.7pt\intop}}

\def\divv{{\rm div}\,}

\def\supp{{\rm supp}\,}

\def\today{${\scriptscriptstyle\number\day-\number\month-\number\year}$}
\footline={{\hfil\rm\the\pageno\hfil${\scriptscriptstyle\rm\jobname}$\ \ \today}}
\def\D{{\Bbb D}}

\def\R{{\Bbb R}}
\def\N{{\Bbb N}}
\def\T{{\Bbb T}}

\def\esssup{{\rm esssup\,}}

\center{\secbf Global regular axially symmetric solutions to the 
Navier-Stokes equations in a periodic cylinder}
\vskip1.5cm

\centerline{\bf Wojciech M. Zaj\c aczkowski}

\vskip1cm
\noindent
Institute of Mathematics, Polish Academy of Sciences,\\
\'Sniadeckich 8, 00-956 Warsaw, Poland\\
E-mail:wz@impan.pl;\\
Institute of Mathematics and Cryptology, Cybernetics Faculty,\\
Military University of Technology, Kaliskiego 2,\\
00-908 Warsaw, Poland
\vskip0.8cm


\vskip1.5cm

\noindent
{\bf Abstract.} 
We examine the axially symmetric solutions to the Navier\--Stokes equations 
in a periodic cylinder with the slip boundary conditions on the lateral 
part of the boundary. Having initial velocity $v(0)\in H^1(\Omega)$ we 
prove the existence of local regular solutions such that 
$v\in W_2^{2,1}(\Omega^{T_*})$, $\nabla p\in L_2(\Omega^{T_*})$, where $T_*$ 
is so small that
$$
cT_*^{1/2}\|v(0)\|_{H^1(\Omega)}\le1.
\leqno(1)
$$
If $v\in W_2^{2,1}(\Omega^{T_*})$ then swirl $u=rv_\varphi$ belongs to 
$C^{1/2,1/4}(\Omega^{T_*})$ and vanishes on the axis of symmetry.

\noindent
Let $v_r$, $v_\varphi$, $v_z$ be the cylindrical components of velocity and 
$\chi=v_{r,z}-v_{z,r}$. Introducing the quantity $(c\ge1)$
$$\eqal{
\alpha&=c\bigg(\|v(0)\|_{H^1(\Omega)}+
\bigg\|{v_\varphi(0)\over\sqrt{r}}\bigg\|_{L_4(\Omega)}+
\bigg\|{v_\varphi(0)\over r}\bigg\|_{L_3(\Omega)}+
\bigg\|{\chi(0)\over r}\bigg\|_{L_2(\Omega)}\bigg)\cr
&\ge\|v(0)\|_{H^1(\Omega)}\cr}
\leqno(2)
$$
we derive the estimate
$$
\|v(T_*)\|_{H^1(\Omega)}\le\|v\|_{V_2^1(\Omega^{T_*})}\le\alpha.
\leqno(3)
$$
Defining $T$ by the relation
$$
cT^{1/2}\alpha\le1
\leqno(4)
$$
we have that $T<T_*$. Then (3) yields that
$$
\|v(T)\|_{H^1(\Omega)}\le\alpha.
\leqno(5)
$$
Starting from time $t=T$ with initial data satisfying (5) and repeating the 
above considerations the existence in the interval $[T,2T]$ is proved and 
also estimate (3) holds, so
$$
\|v(2T)\|_{H^1(\Omega)}\le\alpha.
\leqno(6)
$$
Since estimate (3) is a global a priori estimate the extension can be 
performed step by step in all intervals $(kT,(k+1)T)$, $k\in\N$. Therefore, 
we prove the existence of such solutions that
$$
v\in W_2^{2,1}(\Omega\times(kT,(k+1)T)),\ \nabla p\in L_2(\Omega\times(kT,
(k+1)T)),\ k\in\N\cup\{0\},
$$
where $T$ is determined by (4) with $\alpha$ introduced in (2).

\section{1. Introduction}

We consider the axially symmetric solutions to the problem
$$\eqal{
&v_t+v\cdot\nabla v-\divv\T(v,p)=0\quad &{\rm in}\ \ 
\Omega^T=\Omega\times(0,T),\cr
&\divv v=0\quad &{\rm in}\ \ \Omega^T,\cr
&v\cdot\bar n=0\quad &{\rm on}\ \ S^T=S\times(0,T),\cr
&\bar n\cdot\D(v)\cdot\bar\tau_\alpha=0,\ \ \alpha=1,2,\quad &{\rm on}\ \ 
S_1^T,\cr
&{\rm periodic\ boundary\ conditions}\quad &{\rm on}\ \ S_2^T,\cr
&v|_{t=0}=v_0,\quad &{\rm in}\ \ \Omega,\cr}
\leqno(1.1)
$$
where $\Omega\subset\R^3$ is a cylinder with boundary $S=S_1\cup S_2$, 
$S_1$ is parallel to the axis of the cylinder, but $S_2$ is perpendicular, 
$v=(v_1(x,t),v_2(x,t),\break v_3(x,t))\in\R^3$ 
is the velocity of the fluid, $p=p(x,t)\in\R$ the pressure, 
$x=(x_1,x_2,x_3)$ are the Cartesian 
coordinates such that $x_3$ axis is the axis of the cylinder, $\bar n$ is 
the unit outward normal vector to $S_1$, $\bar\tau_\alpha$, $\alpha=1,2$, is 
the tangent vector to $S_1$ and the dot denotes the scalar product in $\R^3$.

By $\T(v,p)$ we denote the stress tensor of the form
$$
\T(v,p)=\nu\D(v)-pI,
\leqno(1.2)
$$
where $\nu>0$ is the constant viscosity coefficient, $\D(v)$ is the 
dilatation tensor of the form
$$
\D(v)=\{v_{i,x_j}+v_{j,x_i}\}_{ij=1,2}
\leqno(1.3)
$$
and $I$ is the unit matrix.

We are interested to examine axially symmetric solutions to problem (1.1). 
Hence we introduce the cylindrical coordinates $r$, $\varphi$, $z$ by the 
relations $x_1=r\cos\varphi$, $x_2=r\sin\varphi$, $x_3=z$.

\noindent
We assume that $\Omega$ is axially symmetric bounded cylinder with $x_3$ 
axis as the axis of symmetry. Let $R$ and $a$ be given posite numbers. Then
$$
\Omega=\{x\in\R^3:\ r<R,\ |z|<a\}.
$$
Hence
$$
S_1=\{x\in\R^3:\ r=R,\ |z|\le a\}
$$
and
$$
S_2=\{x\in\R^3:\ r\le R,\ z\ {\rm is\ equal\ either}\ -a\ {\rm or}\ a\}.
$$
Let us introduce the vectors
$$
\bar e_r=(\cos\varphi,\sin\varphi,0),\quad
\bar e_\varphi=(-\sin\varphi,\cos\varphi,0),\quad
\bar e_z=(0,0,1),
$$
which are connected with cylindrical coordinates. 
Then the cylindrical coordinates of velocity are defined by the relations
$$
v_r=v\cdot\bar e_r,\quad v_\varphi=v\cdot\bar e_\varphi,\quad 
v_z=v\cdot\bar e_z.
\leqno(1.4)
$$

\Def{1.1.} 
By the axially symmetric solution to problem (1.1) we mean such solution that
$$
v_{r,\varphi}=v_{\varphi,\varphi}=v_{z,\varphi}=p_{,\varphi}=0.
\leqno(1.5)
$$

The aim of this paper is to prove the existence of global regular axially 
symmetric solutions to problem (1.1) with large swirl $u$, which is defined by
$$
u=rv_\varphi.
\leqno(1.6)
$$
We have to recall that behavior of axially symmetric solutions is 
essentially different near and far from the axis of symmetry. An appropriate 
a priori estimate near the axis is found in [Z1]. In this paper we show 
a corresponding estimate outside of the axis. Hence the estimate in the 
whole domain $\Omega$ follows. The estimate, we are looking for, is for the 
norm $\|v(t)\|_{H^1(\Omega)}$, $t\in\R_+$. The estimate is sufficient to prove 
the existence of global regular axially symmetric solutions to (1.1).

\noindent
To formulate the main results of this paper we need

\noindent
{\bf Assumptions:}
\vskip6pt

\item{1.} Let $v_0\in L_2(\Omega)$ and let $d_1$ be a positive constant such 
that\\ 
$\|v_0\|_{L_2(\Omega)}\le d_1$.

\item{2.} Let $u_0\in L_\infty(\Omega)$ and let $d_2$ be a positive constant 
such that\\ 
$\|u_0\|_{L_\infty(\Omega)}\le d_2$.

\item{3.} Let $\zeta_1=\zeta_1(r)$ be a smooth cut off function such that 
$\zeta_1(r)=1$ for $r\le r_0$ and $\zeta_1(r)=0$ for $2r_0\le r<R$. Next 
$\zeta_2(r)=1$ for $r\ge r_0$ and $\zeta_2(r)=0$ for $r\le{r_0\over2}$. 
Finally $\zeta_3(r)=1$ for $r\ge2r_0$ and $\zeta_3(r)=0$ for $r\le r_0$. 
Moreover, $\{\zeta_1(r),\zeta_3(r)\}$ is a partition of unity in $[0,R]$. 
Let $\Omega_{\zeta_i}=\Omega\cap\supp\zeta_i$ and $i=1,2,3$.

\item{4.} Let $u_0\in C^\alpha(\Omega)$, $\alpha\le{1\over2}$. 
Let $r_0$ be so small that
$$
\|u_0\|_{C^\alpha(\Omega_{\zeta_1})}r_0^\alpha\le\root{4}\of{5\over2}\nu
\leqno(1.7)
$$

\item{5.} Assume that
$$
\bigg\|{v_\varphi(0)\over\sqrt{r}}\bigg\|_{L_4(\Omega_{\zeta_1})}<\infty,\quad
\bigg\|{\chi(0)\over r}\bigg\|_{L_2(\Omega_{\zeta_1})}<\infty.
$$

\item{6.} Let
$$\eqal{
&A_1=\varphi(d_1,d_2,1/r_0)\bigg(1+
\bigg\|{v_\varphi(0)\over\sqrt{r}}\bigg\|_{L_4(\Omega_{\zeta_1})}+
\bigg\|{\chi(0)\over r}\bigg\|_{L_2(\Omega_{\zeta_1})}\bigg),\cr
&A_2=\varphi(d_1,1/r_0)(1+\|v_\varphi(0)\|_{L_3(\Omega_{\zeta_2})}+
\|\chi(0)\|_{L_2(\Omega_{\zeta_2})}),\cr
&A=A_1+A_2,\cr
&A_3=\varphi(A_1)(d_1+d_2)+c\bigg(\|v_\varphi(0)\|_{H^1(\Omega)}+
\bigg\|{v_\varphi(0)\over r}\bigg\|_{L_{27\over10}(\Omega)}\bigg),\cr}
$$
where $\varphi$ is an increasing positive function.
\vskip6pt

\proclaim Theorem A. 
Let assumptions 1--6 hold. Let $u_0\in C^\alpha(\Omega)$, $\alpha\le1/2$, 
$v_0\in H^1(\Omega)$. Let $\alpha=A_1+A_2+A_3$. Let $c_*$ be some positive 
constant. Let $T$ be such that $c_*T^{1/2}\alpha\le1$. Then there exists 
a global regular axially symmetric solution to (1.1) such that 
$v\in W_2^{2,1}(\Omega\times(kT,(k+1)T))$, 
$\nabla p\in L_2(\Omega\times(kT,(k+1)T))$  $k\in\N\cup\{0\}$ and
$$\eqal{
&\|v\|_{W_2^{2,1}(\Omega\times(kT,(k+1)T))}+\|\nabla p\|_{L_2(\Omega\times(kT,(k+1)T))}\cr
&\le\varphi(\alpha,\|u_0\|_{C^{1/2}(\Omega)})+c\|v_0\|_{H^1(\Omega)}.\cr}
\leqno(1.9)
$$

\proclaim Theorem B. 
In the above case we have uniqueness of solutions.

Now we shall describe shortly the proof of Theorem A. Having that 
$v(0)\in H^1(\Omega)$ we prove a local existence of solutions to 
problem (1.1) such that $v\in W_2^{2,1}(\Omega^{T_*})$, 
$\nabla p\in L_2(\Omega^{T_*})$, where $T_*$ is so small that
$$
c_*T_*^{1/2}\|v(0)\|_{H^1(\Omega)}\le 1,
\leqno(1.8)
$$
where constant $c_*$ is described in (3.3).

\noindent
Then it follows that $v\in L_r(0,T_*;L_p(\Omega))$, where 
${3\over p}+{2\over r}={1\over2}$. Next Theorem 5.4 from [Z4] implies that 
$u\in C(0,T_*;C^{1/2}(\Omega_{\zeta_1}))$, where $\Omega_{\zeta_1}$ is 
described in the above Assumption 3.

\noindent
Since $u$ is continuous, $\Omega_{\zeta_1}$ belongs to a cylinder with 
radius $2r_0$ and $u$ vanishes on the axis of symmetry we obtain for $r_0$ 
sufficiently small that
$$
\|u\|_{C(0,T_*;C^{1/2}(\Omega_{\zeta_1}))}\le\root{4}\of{5\over4}\nu.
\leqno(1.9)
$$
Then Lemma 6.4 from [Z1] yields the estimate
$$
\|v'\|_{V_2^1(\Omega_{\zeta_1}^{T_*})}\le A_1,
\leqno(1.10)
$$
where $v'=(v_r,v_z)$.

\noindent
Moreover, Lemma 5.2 yields
$$
\|v'\|_{V_2^1(\Omega_{\zeta_3}^{T_*})}\le A_2.
\leqno(1.11)
$$
Since $\{\zeta_1,\zeta_3\}$ is a partition of unity estimates (1.10) and 
(1.11) imply
$$
\|v'\|_{V_2^1(\Omega^{T_*})}\le A_1+A_2.
\leqno(1.12)
$$
The estimate as a priori is valid for any $T_*>0$ (see Lemma 6.4 [Z.1] and 
Lemma 5.2). By Lemma 4.1 we have that
$$
\|v_\varphi(t)\|_{H^1(\Omega)}\le A_3,\quad t\in\R_+.
\leqno(1.13)
$$
From (1.12) and (1.13) we derive
$$
\|v(T_*)\|_{H^1(\Omega)}\le A_1+A_2+A_3\equiv\alpha.
\leqno(1.14)
$$
Let $T$ be such that
$$
c_*T^{1/2}\alpha\le1,
\leqno(1.15)
$$
and $T<T_*$ because $\|v(0)\|_{H^1(\Omega)}\le\alpha$. Hence (1.14) takes the form
$$
\|v(T)\|_{H^1(\Omega)}\le\alpha.
\leqno(1.16)
$$
Using that $\|v(T)\|_{H^1(\Omega)}\le\alpha$ we can repeat the above 
considerations for any interval $[kT,(k+1)T]$, $k\in\N$.

\noindent
The above extension in time of the local solution is possible because $A_i$, 
$i=1,2,3$, do not depend on time and estimates (1.12), (1.13) are global 
a priori estimates (a priori means that it is supposed that there exists 
a solution to (1.1) such that $v\in W_2^{2,1}(\Omega\times\R_+))$.

\section{2. Auxiliary results and notation}
\vskip-12pt

\Def{2.1.} 
We introduce the spaces
$$
V_2^0(\Omega^T)=\{u:\ \|u\|_{L_\infty(0,T;L_2(\Omega))}+
\|\nabla u\|_{L_2(\Omega^T)}<\infty\}
$$
and
$$\eqal{
&V_2^0(\Omega\times(kT,(k+1)T))\cr
&\qquad=\{u:\ \|u\|_{L_\infty(kT,(k+1)T;L_2(\Omega))}+
\|\nabla u\|_{L_2(\Omega\times(kT,(k+1)T))}<\infty\}\cr}
$$
where $k\in\N_0\equiv\N\cup\{0\}$.

\proclaim Lemma 2.2. 
Assume that $v_0\in L_2(\Omega)$ and $d_1>0$ is a constant such that 
$\|v_0\|_{L_2(\Omega)}\le d_1$. Then
$$
\|v(t)\|_{L_2(\Omega)}\le\|v_0\|_{L_2(\Omega)}\le d_1.
\leqno(2.1)
$$
Assume that $d_0>0$ is a constant such that 
$\left|\intop_\Omega u_0dx\right|\le d_0$. Then
$$
\bigg|\intop_\Omega u(t)dx\bigg|=\bigg|\intop_\Omega u_0dx\bigg|\le d_0.
\leqno(2.2)
$$
Assume that $T>0$ is fixed. Then
$$
\|v\|_{V_2^0(\Omega\times(kT,t))}^2\le{2\over\min\{1,\nu_2\}}
\bigg({d_0^2\over\nu_1}e^{\nu_1T}+\|v(kT)\|_{L_2(\Omega)}^2\bigg),
\leqno(2.3)
$$
where $t\in[kT,(k+1)T]$, $k\in\N_0$, $\nu_*=\nu_1+\nu_2$, $\nu_i>0$, $i=1,2$, 
$\nu_*={\nu\over c_k}$, $c_k$ is the constant from the Korn 
inequality (2.5).

\Proof 
Multiplying $(1.1)_1$ by $v$, integrating over $\Omega$, using $(1.1)_2$ and 
boundary conditions we obtain
$$
{d\over dt}\|v\|_{L_2(\Omega)}^2+\nu E_\Omega(v)=0,
\leqno(2.4)
$$
where
$$
E_\Omega(v)=\|\D(v)\|_{L_2(\Omega)}^2.
$$
Integrating (2.4) with respect to time and omitting the second expression on 
the l.h.s. implies (2.1).

\noindent
From [Z3, Ch. 4, Lemma 2.4] we have the following Korn inequality
$$
\|v\|_{H^1(\Omega)}^2\le c_k\bigg(E_\Omega(v)+
\bigg|\intop_\Omega v\cdot\eta dx\bigg|^2\bigg),
\leqno(2.5)
$$
where $\eta=(-x_2,x_1,0)$, $v\cdot\eta=rv_\varphi=u$.

Now we calculate the last term on the r.h.s. of (2.5). Multiplying $(1.1)_1$ 
by $\eta$, integrating over $\Omega$, using $(1.1)_2$ we obtain
$$
{d\over dt}\intop_\Omega v\cdot\eta dx-\intop_\Omega v_iv_j\nabla_i\eta_jdx+
\intop_\Omega\T_{ij}\nabla_i\eta_jdx=0,
\leqno(2.6)
$$
where the summation convention over the repeated indices is assumed.

\noindent
Using that $\nabla\eta$ is antisymmetric tensor we see that (2.6) implies
$$
{d\over dt}\intop_\Omega v\cdot\eta dx=0
\leqno(2.7)
$$
Hence
$$
\intop_\Omega u(t)dx=\intop_\Omega u_0dx.
\leqno(2.8)
$$
Employing (2.5) and (2.8) in (2.4) yields
$$
{d\over dt}\|v\|_{L_2(\Omega)}^2+\nu_*\|v\|_{H^1(\Omega)}^2\le d_0^2,
\leqno(2.9)
$$
where the inequality (2.2) was used and $\nu_*={\nu\over c_k}$.

Let $\nu_*=\nu_1+\nu_2$, $\nu_i>0$, $i=1,2$. Then (2.9) implies
$$
{d\over dt}(\|v\|_{L_2(\Omega)}^2e^{\nu_1t})+\nu_2
\|v\|_{H^1(\Omega)}^2e^{\nu_1t}\le d_0^2e^{\nu_1t}.
\leqno(2.10)
$$
Integrating (2.10) with respect to time from $kT$ to $t$ yields
$$\eqal{
&\|v(t)\|_{L_2(\Omega)}^2+\nu_2e^{-\nu_1t}\intop_{kT}^t
\|v(t')\|_{H^1(\Omega)}^2e^{\nu_1t'}dt'\cr
&\le{d_0^2\over\nu_1}+e^{-\nu_1(t-kT)}\|v(kT)\|_{L_2(\Omega)}^2,\cr}
\leqno(2.11)
$$
where $t\in[kT,(k+1)T]$.

Omitting the first norm on the l.h.s. of (2.11) we obtain
$$\eqal{
&\nu_2\intop_{kT}^t\|v(t')\|_{H^1(\Omega)}^2dt'\le{d_0^2\over\nu_1}
e^{\nu_1(t-kT)}+\|v(kT)\|_{L_2(\Omega)}^2,\cr
&t\in[kT,(k+1)T].\cr}
\leqno(2.12)
$$
Hence, (2.11) and (2.12) imply (2.3). This concludes the proof.

\noindent
From properties of $V_2^0(\Omega^T)$ we have
$$
\|v\|_{L_q(0,T;L_p(\Omega))}\le c_2\|v\|_{V_2^0(\Omega^T)},
\leqno(2.13)
$$
where
$$
{3\over p}+{2\over q}\ge{3\over2}.
$$
Let us consider the problem
$$\eqal{
&v_t-\divv\T(v,p)=f\quad &{\rm in}\ \ \Omega^T,\cr
&\divv v=0\quad &{\rm in}\ \ \Omega^T,\cr
&v\cdot\bar n=0,\ \ \bar n\cdot\D(v)\cdot\bar\tau_\alpha=0,\ \ \alpha=1,2\quad
&{\rm on}\ \ S^T,\cr
&v|_{t=0}=v_0\quad &{\rm in}\ \ \Omega.\cr}
\leqno(2.14)
$$
From [ZZ] we have

\proclaim Lemma 2.3.
Assume that $f\in L_s(\Omega^T)$, $v_0\in W_s^{2-2/s}(\Omega)$, 
$s\in(1,\infty)$, $S_1\subset C^2$. Then there exists a solution to problem 
(2.14) such that $v\in W_s^{2,1}(\Omega^T)$, $\nabla p\in L_s(\Omega^T)$ 
and there exists a constant $c=c(\Omega,s)$ such that
$$
\|v\|_{W_s^{2,1}(\Omega^T)}+\|\nabla p\|_{L_s(\Omega^T)}\le c(\Omega,s)
(\|f\|_{L_s(\Omega^T)}+\|v_0\|_{W_s^{2-2/s}(\Omega)}).
\leqno(2.15)
$$

\noindent
In view of Definition 1.1 equations $(1.1)_{1,2}$ for the axially symmetric 
solutions assume the form
$$
v_{r,t}+v\cdot\nabla v_r-{v_\varphi^2\over r}-\nu\Delta v_r+\nu
{v_r\over r^2}=-p_{,r},
\leqno(2.16)
$$
$$
v_{\varphi,t}+v\cdot\nabla v_\varphi+{v_r\over r}v_\varphi-\nu\Delta v_\varphi+
\nu{v_\varphi\over r^2}=0
\leqno(2.17)
$$
$$
v_{z,t}+v\cdot\nabla v_z-\nu\Delta v_z=-p_{,z},
\leqno(2.18)
$$
$$
v_{r,r}+v_{z,z}=-{v_r\over r},
\leqno(2.19)
$$
where $v\cdot\nabla=v_r\partial_r+v_z\partial_z$, 
$\Delta u={1\over r}(ru_{,r})_{,r}+u_{,zz}$.

\noindent
The slip-boundary conditions $(1.1)_{3,4}$ on $S_1$ imply (see [Z3, Ch. 2])
$$
v_r=0,\quad v_{z,r}=0,\quad v_{\varphi,r}={1\over R}v_\varphi\quad 
{\rm on}\ \ S_1.
\leqno(2.20)
$$

\proclaim Lemma 2.4. 
For the weak solutions to problem (1.1) we have the estimate
$$\eqal{
&\|v\|_{V_2^0(\Omega\times(kT,t))}^2+
\bigg\|{v_r\over r}\bigg\|_{L_2(\Omega\times(kT,t))}^2+
\bigg\|{v_\varphi\over r}\bigg\|_{L_2(\Omega\times(kT,t))}^2\cr
&\le c_0\|v_0\|_{L_2(\Omega)}^2\equiv d_1^2,\cr}
\leqno(2.21)
$$
where $c_0=c(T+1)$, $t\in[kT,(k+1)T]$, $k\in\N_0$ and $c$ is the constant 
from imbedding $H^1(\Omega)\subset L_2(S_1)$.

\Proof 
Multiplying (2.17) by $v_\varphi$, integrating over $\Omega$ and using 
boundary conditions (2.20) yields
$$\eqal{
&{1\over2}{d\over dt}\intop_\Omega v_\varphi^2dx+\intop_\Omega
{v_r\over r}v_\varphi^2dx+\nu\intop_\Omega(v_{\varphi,r}^2+
v_{\varphi,z}^2)dx-\nu\intop_{-a}^av_\varphi^2dz\cr
&\quad+\nu\intop_\Omega{v_\varphi^2\over r^2}dx=0.\cr}
$$
Multiplying (2.16) by $v_r$, integrating over $\Omega$ and using boundary 
conditions (2.20) implies
$$\eqal{
&{1\over2}{d\over dt}\intop_\Omega v_r^2dx-\intop_\Omega
{v_\varphi^2\over r}v_rdx+\nu\intop_\Omega(v_{r,r}^2+v_{r,z}^2)dx+
\nu\intop_\Omega{v_r^2\over r^2}dx\cr
&=-\intop_\Omega p_{,r}v_rdx.\cr}
$$
Multiplying (2.18) by $v_z$, integrating over $\Omega$ and using boundary 
conditions (2.20) we obtain
$$
{1\over2}{d\over dt}\intop_\Omega v_z^2dx+\nu\intop_\Omega(v_{z,r}^2+v_{z,z}^2)
dx=-\intop_\Omega p_{,z}v_zdx.
$$
Adding the above equations, using (2.19) and the inequality 
$$
\|v_\varphi\|_{L_2(S_1)}^2\le\varepsilon\|\nabla v_\varphi\|_{L_2(\Omega)}^2+
c(1/\varepsilon)\|v_\varphi\|_{L_2(\Omega)}^2
$$
we obtain for sufficiently small $\varepsilon$,
$$\eqal{
&{1\over2}{d\over dt}\intop_\Omega(v_r^2+v_\varphi^2+v_z^2)dx+{\nu\over2}
\intop_\Omega(v_{r,r}^2+v_{r,z}^2+v_{\varphi,r}^2+v_{\varphi,z}^2+v_{z,r}^2+
v_{z,z}^2)dx\cr
&\quad+\nu\intop_\Omega\bigg({v_r^2\over r^2}+{v_\varphi^2\over r^2}\bigg)dx
\le c\intop_\Omega v_\varphi^2dx.\cr}
\leqno(2.22)
$$
Integrating (2.22) with respect to time from $kT$ to $t\in(kT,(k+1)T]$, 
$k\in\N_0$ and using (2.1) yields
$$\eqal{
&\|v(t)\|_{L_2(\Omega)}^2+\nu\intop_{kT}^t\intop_\Omega(|v_{,r}|^2+|v_{,z}|^2)
dxdt'+\nu\intop_{kT}^t\intop_\Omega\bigg({v_r^2\over r^2}+
{v_z^2\over r^2}\bigg)dxdt'\cr
&\le cTd_1^2+\|v(kT)\|_{L_2(\Omega)}^2.\cr}
$$
Applying again (2.1) gives (2.21). This ends the proof.
\kwadrat

\noindent
Let us introduce the component of vorticity
$$
\chi=v_{r,z}-v_{z,r}.
\leqno(2.23)
$$
Then $\chi$ is a solution to the problem (see [Z1, Z3])
$$\eqal{
&\chi_t+v\cdot\nabla\chi-{v_r\over r}\chi-\nu\bigg[\bigg(r
\bigg({\chi\over r}\bigg)_{,r}\bigg)_{,r}+\chi_{,zz}+2
\bigg({\chi\over r}\bigg)_{,r}\bigg]\cr
&=2{v_\varphi v_{\varphi,z}\over r},\cr
&\chi|_{S_1}=0,\ \ \chi|_{S_2}\ \ {\rm satisfies\ periodic\ boundary\ 
conditions,}\cr
&\chi|_{t=0}=\chi_0.\cr}
\leqno(2.24)
$$

\Def{2.5.} 
By $V_2^k(\Omega^T)$, $k\in\N\cup\{0\}$, we denote the space of functions 
with the following finite norm
$$
\|u\|_{V_2^k(\Omega^T)}=\esssup_{t\in[0,T]}\|u\|_{H^k(\Omega)}+
\|\nabla u\|_{L_2(0,T;H^k(\Omega))},
$$
where
$$
\|u\|_{H^k(\Omega)}=\bigg(\sum_{|\alpha|\le k}\intop_\Omega
|D_x^\alpha u|^2dx\bigg)^{1/2},
$$
$D_x^\alpha=\partial_{x_1}^{\alpha_1}\partial_{x_2}^{\alpha_2}
\partial_{x_3}^{\alpha_3}$, $\alpha=\alpha_1+\alpha_2+\alpha_3$, 
$\alpha_i\in\N\cup\{0\}$, $i=1,2,3$.
\goodbreak

\Remark{2.6.} 
From (2.21) we have that swirl $u=rv_\varphi$ satisfies
$$
\bigg\|{u\over r^2}\bigg\|_{L_2(\Omega^T)}\le d_1,
\leqno(2.25)
$$
but from [Z4] it follows that $u\in C^{\alpha,\alpha/2}(\Omega^T)$, 
$\alpha>0$. Hence $u|_{r=0}=0$.

\proclaim Lemma 2.7. 
Let $d_1$, $d_2$ be given positive constants. Let 
$\|v_0\|_{L_2(\Omega)}\le d_1$, $\|u\|_{L_\infty(\Omega^T)}\le d_2$ (see 
Assumptions 1, 2). Then
$$
\|v_\varphi\|_{L_4(\Omega^T)}\le d_1^{1/2}d_2^{1/2}.
\leqno(2.26)
$$
Proof follows from application of Lemma 2.4 and Lemma 2.1 from [Z4].

\section{3. Local existence}

We prove a local existence of solutions to problem (1.1) by the 
Leray\--Schauder fixed point theorem applying the old idea of 
O. A. Ladyzhenskaya (see [L1, Ch. 4, Theorem 1']). 

\noindent
Let us consider the problem
$$\eqal{
&v_t-\nu\Delta v+\nabla p=-v\cdot\nabla v\quad &{\rm in}\ \ 
\Omega_\varepsilon^T,\cr
&\divv v=0\quad &{\rm in}\ \ \Omega_\varepsilon^T,\cr
&\bar n\cdot v=0,\ \ \bar n\cdot\T(v,p)\cdot\bar\tau_\alpha=0,\ \ 
\alpha=1,2,\quad &{\rm on}\ \ S_1^T\cup S_\varepsilon^T,\cr
&{\rm periodic\ boundary\ conditions}\quad &{\rm on}\ \ S_2^T,\cr
&v|_{t=0}=v(0)\quad &{\rm in}\ \ \Omega_\varepsilon,\cr}
\leqno(3.1)
$$
where $\Omega_\varepsilon=\{x\in\Omega\colon r>\varepsilon\}$, 
$S_\varepsilon=\{x\in\R^3\colon r=\varepsilon,\ |z|<a\}$. Setting 
$\varepsilon=0$ we obtain problem (1.1).

\noindent
In view of [S, ZZ] there exists a constant $c_0$ such that solutions to (3.1) 
satisfy the inequality
$$
\|v\|_{W_2^{2,1}(\Omega_\varepsilon^T)}\le c_0(\|v\cdot\nabla v\|_{L_2(\Omega_\varepsilon^T)}+
\|v(0)\|_{H^1(\Omega_\varepsilon)}).
\leqno(3.2)
$$

\proclaim Theorem 3.1. 
Assume that $v(0)\in H^1(\Omega)$. Assume that $T$ so small that
$$
c_*T^{1/2}\|v(0)\|_{H^1(\Omega)}\le 1
\leqno(3.3)
$$
where $c_*$ is such constant that (3.10) is satisfied. Then there exists 
a solution to problem (1.1) such that $v\in W_2^{2,1}(\Omega^T)$, 
$\nabla p\in L_2(\Omega^T)$, where $T$ satisfies (3.3) and
$$
\|v\|_{W_2^{2,1}(\Omega^T)}+\|\nabla p\|_{L_2(\Omega^T)}\le 8c_0
\|v(0)\|_{H^1(\Omega)},
\leqno(3.4)
$$
where $c_0$ is taken from (3.2).

\Proof 
Since we a going to prove the existence of solutions to problem (3.1) by the 
Leray-Schauder fixed point theorem we restrict the proof to show only estimate 
(3.4) because other steps of it are clear.
We shall skip the index $\varepsilon$ in $\Omega_\varepsilon$ for simplicity

Now we examine the first term on the r.h.s. of (3.2). We estimate it by
$$\eqal{
&\bigg(\intop_0^Tdt\intop_\Omega|v\cdot\nabla v|^2dx\bigg)^{1/2}\le
\bigg(\intop_0^T\|v(t)\|_{L_\infty(\Omega)}^2
\|\nabla v(t)\|_{L_2(\Omega)}^2dt\bigg)^{1/2}\cr
&\le\sup_t\|\nabla v(t)\|_{L_2(\Omega)}\bigg(\intop_0^T
\|v(t)\|_{L_\infty(\Omega)}^2dt\bigg)^{1/2}\equiv I_1.\cr}
$$
From [BIN, Ch. 3, Sect. 15] we have
$$
\|v\|_{L_\infty(\Omega)}\le c_1\|v_{xx}\|_{L_2(\Omega)}^{3/4}
\|v\|_{L_2(\Omega)}^{1/4}+c_1\|v\|_{L_2(\Omega)}.
\leqno(3.5)
$$
Using (3.5) in $I_1$ yields
$$\eqal{
I_1&\le\sup_t\|\nabla v\|_{L_2(\Omega)}(c_1T^{1/8}
\|v_{xx}\|_{L_2(\Omega^T)}^{3/4}\sup_t\|v(t)\|_{L_2(\Omega)}^{1/4}\cr
&\quad+c_1T^{1/2}\sup_t\|v(t)\|_{L_2(\Omega)})\equiv I_2.\cr}
$$
Employing
$$
\sup_t\|\nabla v\|_{L_2(\Omega)}\le\sup_t\|v(t)\|_{H^1(\Omega)}\le c_2
(\|v\|_{W_2^{2,1}(\Omega^T)}+\|v(0)\|_{H^1(\Omega)})
$$
and the energy estimate (2.1) in $I_2$, we obtain from (3.2) the inequality
$$\eqal{
&\|v\|_{W_2^{2,1}(\Omega^T)}+\|\nabla p\|_{L_2(\Omega^T)}\le c_0c_2
(\|v\|_{W_2^{2,1}(\Omega^T)}+\|v(0)\|_{H^1(\Omega)})\cdot\cr
&\quad\cdot(c_1T^{1/8}\|v_{xx}\|_{L_2(\Omega^T)}^{3/4}d_1^{1/4}+c_1T^{1/2}d_1)+
c_0\|v(0)\|_{H^1(\Omega)}.\cr}
$$
Assuming that $T$ is so small that
$$
c_0c_2c_1(T^{1/8}\|v_{xx}\|_{L_2(\Omega^T)}^{3/4}d_1^{1/4}+T^{1/2}d_1)\le1/2
\leqno(3.6)
$$
we obtain the inequality
$$\eqal{
&\|v\|_{W_2^{2,1}(\Omega^T)}+\|\nabla p\|_{L_2(\Omega^T)}\le2c_0c_1c_2
(T^{1/2}\|v\|_{W_2^{2,1}(\Omega^T)}^{3/4}d_1^{1/4}\cr
&\quad+T^{1/2}d_1)\|v(0)\|_{H^1(\Omega)}+2c_0\|v(0)\|_{H^1(\Omega)}.\cr}
\leqno(3.7)
$$
By the Young inequality applied to the first term on the r.h.s. of (3.7) 
we have
$$\eqal{
&\|v\|_{W_2^{2,1}(\Omega^T)}+\|\nabla p\|_{L_2(\Omega^T)}\le
{\varepsilon^{4/3}\over4/3}\|v\|_{W_2^{2,1}(\Omega^T)}\cr
&\quad+{1\over4\varepsilon^4}(2c_0c_1c_2T^{1/2}d_1^{1/4}
\|v(0)\|_{H^1(\Omega)})^4\cr
&\quad+2c_0(c_1c_2T^{1/2}d_1+1)\|v(0)\|_{H^1(\Omega)}.\cr}
$$
Setting ${\varepsilon^{4/3}\over4/3}={1\over2}$ so 
$\varepsilon=\big({2\over3}\big)^{3/4}$ the above inequality yields
$$\eqal{
&{1\over2}\|v\|_{W_2^{2,1}(\Omega^T)}+\|\nabla p\|_{L_2(\Omega^T)}\le
{27\over32}(2c_0c_1c_2)^4T^2d_1\|v(0)\|_{H^1(\Omega)}^4\cr
&\quad+2c_0(c_1c_2T^{1/2}d_1+1)\|v(0)\|_{H^1(\Omega)}.\cr}
$$
Simplifying the inequality yields
$$\eqal{
&\|v\|_{W_2^{2,1}(\Omega^T)}+\|\nabla p\|_{L_2(\Omega^T)}\le[27(c_0c_1c_2)^4
T^2d_1\|v(0)\|_{H^1(\Omega)}^3\cr
&\quad+4c_0c_1c_2T^{1/2}d_1]\|v(0)\|_{H^1(\Omega)}+4c_0\|v(0)\|_{H^1(\Omega)}.
\cr}
$$
Assuming that $T$ is so small that
$$
27(c_0c_1c_2)^4T^2d_1\|v(0)\|_{H^1(\Omega)}^3+4c_0c_1c_2T^{1/2}d_1\le4c_0
\leqno(3.8)
$$
we obtain (3.4) after passing with $\varepsilon$ to 0.

\noindent
Using (3.4) in (3.6) we assume the stronger restriction
$$
c_0c_1c_2[(4c_0)^{3/4}T^{1/8}d_1^{1/4}+T^{1/2}d_1]\le1/2.
\leqno(3.9)
$$
Using that $d_1\le\|v(0)\|_{H^1(\Omega)}$ wee find constants $\bar c_i$, 
$i=1,\dots,4$, such that (3.8) and (3.9) can be expressed in the form
$$\eqal{
&\bar c_1T^2\|v(0)\|_{H^1(\Omega)}^4+\bar c_2T^{1/2}
\|v(0)\|_{H^1(\Omega)}\le{1\over2},\cr
&\bar c_3T^{1/8}\|v(0)\|_{H^1(\Omega)}^{1/4}+\bar c_4T^{1/2}
\|v(0)\|_{H^1(\Omega)}\le{1\over2}.\cr}
\leqno(3.10)
$$
Therefore there exists a constant $c_*$ satisfying
$$
c_*T^{1/2}\|v(0)\|_{H^1(\Omega)}\le1
$$
and such that (3.10) is satisfied. This concludes the proof.

\section{4. Estimate for the angular component of velocity}

In this section we prove the following estimate
$$
\|v_\varphi(t)\|_{H^1(\Omega)}\le d_5,
\leqno(4.1)
$$
where $d_5$ does not depend on $t$.

\noindent
For this purpose we consider the problem for $v_\varphi$
$$\eqal{
&v_{\varphi,t}-\nu\Delta v_\varphi+v'\cdot\nabla v_\varphi+{v_r\over r}
v_\varphi+\nu{v_\varphi\over r^2}=0\quad &{\rm in}\ \ \Omega^T,\cr
&v_{\varphi,r}={1\over r}v_\varphi\quad &{\rm on}\ \ S_1^T,\cr
&v_\varphi|_{S_2}\ -\ {\rm is\ periodic}\cr
&v_\varphi|_{t=0}=v_\varphi(0)\quad &{\rm in}\ \ \Omega.\cr}
\leqno(4.2)
$$

\proclaim Lemma 4.1. 
Assume that there exists a local solution described by Theorem 3.1. Assume 
also that $A_0=A_1+A_2$ is finite (see assumptions from Section 1). 
Assume that
$$\eqal{
d_5&=\varphi(A_0)(d_1+d_2)+c(\|v_\varphi(0)\|_{H^1(\Omega)}+
\|v_\varphi(0)\|_{W_{5/2}^{1/5}(\Omega)}\cr
&\quad+\bigg\|{v_\varphi(0)\over r}\bigg\|_{L_{27\over10}(\Omega)}\le A_3\cr}
$$
is finite and does not depend on $t$. Then (4.1) holds for any $t\in\R_+$.
\goodbreak

\Proof 
Multiplying $(4.2)_1$ by $v_{\varphi,t}$ and integrating the result over 
$\Omega$ yields
$$\eqal{
&\intop_\Omega v_{\varphi,t}^2dx+{\nu\over2}{d\over dt}\intop_\Omega
|\nabla v_\varphi|^2dx+{\nu\over2}{d\over dt}\intop_\Omega
{v_\varphi^2\over r^2}dx\cr
&\quad-{\nu\over2}{d\over dt}\intop_{-a}^av_\varphi^2|_{r=R}dz\le
{\varepsilon_1\over2}\intop_\Omega v_{\varphi,t}^2dx+{1\over2\varepsilon_1}
\intop_\Omega|v'\cdot\nabla v_\varphi|^2dx\cr
&\quad+{\varepsilon_2\over2}\intop_\Omega v_{\varphi,t}^2dx+
{1\over2\varepsilon_2}\intop_\Omega\bigg|{v_r\over r}\bigg|^2
v_\varphi^2dx.\cr}
\leqno(4.3)
$$
Setting $\varepsilon_1=\varepsilon_2={1\over2}$ we get
$$\eqal{
&{1\over2}\intop_\Omega v_{\varphi,t}^2dx+{\nu\over2}{d\over dt}\intop_\Omega
|\nabla v_\varphi|^2dx+{\nu\over2}{d\over dt}\intop_\Omega
{v_\varphi^2\over r^2}dx\cr
&\le{\nu\over2}{d\over dt}\intop_{-a}^av_\varphi^2|_{r=R}dz+\intop_\Omega
|v'\cdot\nabla v_\varphi|^2dx+\intop\bigg|{v_r\over r}\bigg|^2
v_\varphi^2dx.\cr}
\leqno(4.4)
$$
Integrating (4.4) with respect to time and using estimates
$$
\|v'\|_{L_{10}(\Omega^T)}\le\varphi(A_0)
\leqno(4.5)
$$
from [Z1, (6.19) and (6.20)],
$$
\bigg\|{v_r\over r}\bigg\|_{L_{10}(\Omega^T)}\le\varphi(A_0)
\leqno(4.6)
$$
from [Z1, (6.38)], and (see (4.17))
$$
\|v_\varphi|_{r=R}\|_{L_\infty(\Omega^T)}\le c\|u_0\|_{L_\infty(\Omega)}
\leqno(4.7)
$$
we obtain
$$\eqal{
&\intop_{\Omega^t}v_{\varphi,t}^2dxdt+\nu\intop_\Omega|\nabla v_\varphi(t)|^2dx
+\nu\intop_\Omega{v_\varphi^2(t)\over r^2}dx\cr
&\le\nu\intop_{-a}^av_\varphi^2|_{r=R}dz+\nu\intop_\Omega
|\nabla v_\varphi(0)|^2dx+\nu\intop_\Omega{v_\varphi^2(0)\over r^2}dx\cr
&\quad+\varphi(A_0)(\|\nabla v_\varphi\|_{L_{5/2}(\Omega^t)}^2+
\|v_\varphi\|_{L_{5/2}(\Omega^t)}^2),\cr}
\leqno(4.8)
$$
where the first integral on the r.h.s. is bounded by
$$
c\|u(0)\|_{L_\infty(\Omega)}.
$$
To estimate the norms from the last term on the r.h.s. of (4.8) we introduce 
the Green function to problem (4.2). Let us denote it by $G$. Then we can 
express (4.2) in the integral form
$$\eqal{
v_\varphi(x,t)&=\intop_{\Omega^t}\nabla_yG(x-y,t-\tau)v'v_\varphi dyd\tau\cr
&\quad-\intop_{\Omega^t}G(x-y,t-\tau)\bigg({v_r\over r}v_\varphi+\nu
{v_\varphi\over r^2}\bigg)dyd\tau\cr
&\quad+\intop_\Omega G(x-y,t)v_\varphi(y,0)dy+\intop_{S_1^t}G(x-z,t-\tau)
v_\varphi(R,z)dzd\tau.\cr}
\leqno(4.9)
$$
From (4.9) we have
$$\eqal{
&\|v_\varphi\|_{W_\sigma^{1,1/2}(\Omega^T)}\le c
\bigg(\|v'v_\varphi\|_{L_\sigma(\Omega^T)}+
\bigg\|{v_r\over r}v_\varphi\bigg\|_{L_{5\sigma\over5+\sigma}(\Omega^T)}\cr
&\quad+\bigg\|{v_\varphi\over r^2}\bigg\|_{L_{5\sigma\over5+\sigma}(\Omega^T)}+
\|v_\varphi(0)\|_{W_\sigma^{1-{2\over\sigma}}(\Omega)}+
c\|u(0)\|_{L_\infty(\Omega)}\bigg),\cr}
\leqno(4.10)
$$
where the last norm estimates the last integral on the r.h.s. of (4.9).

\noindent
To estimate the norm $\|\nabla v_\varphi\|_{L_{5/2}(\Omega^T)}$ appearing 
on the r.h.s. of (4.8) we assume $\sigma={5\over2}$.

\noindent
Now we estimate the particular terms on the r.h.s. of (4.10). We bound the 
first term by
$$\eqal{
&\bigg\|{v'\over r}rv_\varphi\bigg\|_{L_{5/2}(\Omega^T)}=
\bigg\|{v'\over r}u\bigg\|_{L_{5/2}(\Omega^T)}\le\|u(0)\|_{L_\infty(\Omega)}
\bigg\|{v'\over r}\bigg\|_{L_{5/2}(\Omega^T)}\cr
&\le\varphi(A_0)\|u(0)\|_{L_\infty(\Omega)},\cr}
$$
where we used (6.19), (6.20), (6.38) from [Z1] (see (4.6), (4.7)).

\noindent
The second integral on the r.h.s. of (4.10) is estimated by
$$
\bigg\|v_r{v_\varphi\over r}\bigg\|_{L_{5\over3}(\Omega^T)}\le
\|v_r\|_{L_{10}(\Omega^T)}
\bigg\|{v_\varphi\over r}\bigg\|_{L_2(\Omega^T)}\le\varphi(A_0)d_1
$$
where we used (6.19) and (6.20) from [Z1].

\noindent
Assuming that $v_\varphi\ge1$, because otherwise we have regularity of axially 
symmetric solutions (see [NP1, NP2, P]), the third integral on the r.h.s. of 
(4.10) is bounded by
$$
\bigg\|{v_\varphi\over r}\bigg\|_{L_{10\over3}(\Omega^T)}.
\leqno(4.11)
$$
Finally, the fourth term on the r.h.s. of (4.10) equals
$$
\|v_\varphi(0)\|_{W_{5\over2}^{1-4/5}(\Omega)}
$$
To estimate (4.11) we introduce the quantity $\omega={v_\varphi\over r}$ and 
multiply $(4.2)_1$ by ${1\over r}$. Then we obtain the following equation 
for $\omega$
$$\eqal{
&\omega_{,t}+v'\cdot\nabla\omega+{2v_r\over r}\omega-\nu\Delta\omega-
{2\nu\over r}\omega_{,r}=0\quad &{\rm in}\ \ \Omega^T,\cr
&\omega_{,r}=0\quad &{\rm on}\ \ S_1^T,\cr
&\omega|_{S_2}\ -\ {\rm periodic\ with\ respect\ to}\ z,\cr
&\omega|_{t=0}=\omega(0)\quad &{\rm in}\ \ \Omega.\cr}
\leqno(4.12)
$$
Multiplying (4.12) by $\omega|\omega|^{s-2}$ and integrating over $\Omega$ 
yields
$$\eqal{
&{1\over s}{d\over dt}\intop_\Omega|\omega|^sdx+{4\nu(s-1)\over s^2}
\intop_\Omega|\nabla|\omega^{s/2}|^2dx+{2\nu\over s}\intop_{-a}^a
|\omega|^s|_{r=0}dz\cr
&=-2\intop_\Omega{v_r\over r}|\omega|^sdx+{2\nu\over s}\intop_{-a}^a
|\omega|_{r=R}^sdz.\cr}
\leqno(4.13)
$$
Applying the trace estimate
$$
\intop_{-a}^a|\,|\omega^{s/2}|^2dz\le\varepsilon\intop_\Omega
|\nabla|\omega|^{s/2}|^2dx+{c\over\varepsilon}\intop_\Omega|\omega|^sdx
$$
to the last term on the r.h.s. of (4.13) and assuming that 
$\varepsilon={s-1\over s}$ we obtain
$$\eqal{
&{1\over s}{d\over dt}\intop_\Omega|\omega|^sdx+{2\nu(s-1)\over s^2}
\intop_\Omega|\nabla|\omega|^{s/2}|^2dx\cr
&\le-2\intop_\Omega{v_r\over r}|\omega|^sdx+{cs\over s-1}\intop_\Omega
|\omega|^sdx.\cr}
\leqno(4.14)
$$
Integrating (4.14) with respect to time we get
$$\eqal{
&{1\over s}\intop_\Omega|\omega|^sdx+{2\nu(s-1)\over s^2}\intop_{\Omega^t}
|\nabla|\omega|^{s/2}|^2dxdt'\cr
&\le2\intop_{\Omega^t}\bigg|{v_r\over r}\bigg||\omega|^sdxdt'+
{2\nu\over s}\intop_{\Omega^t}|\omega|^sdxdt'+{1\over s}\intop_\Omega
|\omega(0)|^sdx,\cr}
\leqno(4.15)
$$
where the first integral on the r.h.s. is bounded by
$$
2\bigg\|{v_r\over r}\bigg\|_{L_{10}(\Omega^t)}
\|\omega\|_{L_{{10\over9}s}(\Omega^t)}^s\le\varphi(A_0)
\|\omega\|_{L_{{10\over9}s}(\Omega^t)}^s,
$$
where (6.38) from [Z1] was used.

\noindent
Setting $s={9\over5}$ and using the energy estimate (2.21) we derive
$$
\|\omega\|_{L_3(\Omega^t)}\le\varphi(A_0,d_1)d_1+c\|\omega(0)\|_{L_{9\over5}(\Omega)},
\leqno(4.16)
$$
where $\varphi$ and $c$ do not depend on $t$.\\
Next, setting $s={27\over10}$ in (4.15) and employing (4.16) we obtain from 
(4.15) the estimate
$$
\|\omega\|_{L_{9/2}(\Omega^t)}\le\varphi(A_0,d_1)d_1+c
\|\omega(0)\|_{L_{27\over10}(\Omega)},
\leqno(4.17)
$$
where $\varphi$ and $c$ do not depend on $t$.\\
The inequality is sufficient to estimate (4.11). Summarizing, we have
$$\eqal{
&\|v_\varphi\|_{L_{5/2}(\Omega^T)}+\|\nabla v_\varphi\|_{L_{5/2}(\Omega^T)}\le
\varphi(A_0,d_1)[\|u(0)\|_{L_\infty(\Omega)}+d_1]\cr
&\quad+c\bigg(\|v_\varphi(0)\|_{W_{5/2}^{1/5}(\Omega)}+
\bigg\|{v_\varphi(0)\over r}\bigg\|_{L_{27\over10}(\Omega)}+
\|u(0)\|_{L_\infty(\Omega)}\bigg).\cr}
\leqno(4.18)
$$
Using (4.18) in (4.8), applying the trace estimate to the first integral on 
the r.h.s. of (4.8) and employing (2.21) we obtain
$$\eqal{
&\|v_\varphi(t)\|_{H^1(\Omega)}\le\varphi(A_0,d_1)(\|u(0)\|_{L_\infty(\Omega)}+d_1)\cr
&\quad+c\bigg(\|v_\varphi(0)\|_{H^1(\Omega)}+
\|v_\varphi(0)\|_{W_{5/2}^{1/5}(\Omega)}+
\bigg\|{v_\varphi(0)\over r}\bigg\|_{L_{27/10}(\Omega)}\cr
&\quad+\|u(0)\|_{L_\infty(\Omega)}\bigg).\cr}
\leqno(4.19)
$$
This concludes the proof.

\section{5. Existence}

Let $\zeta_1=\zeta_1(r)$ be a smooth function such that $\zeta_1(r)=0$ 
for $r\ge2r_0$ and $\zeta_1(r)=1$ for $r\le r_0$. By $\zeta_2=\zeta_2(r)$ 
we denote such smooth function that $\zeta_2(r)=1$ for $r\ge r_0$ and 
$\zeta_2(r)=0$ for $r\le r_0/2$. Let us introduce the notation:
$$
v^{'(i)}=v'\zeta_i,\quad \chi^{(i)}=\chi\zeta_i^2,\quad
v_\varphi^{(k)}=v_\varphi\zeta_i,\ \ i=1,2.
$$
From Lemma 6.2 in [Z1] we have

\proclaim Lemma 5.1. 
Assume that $v$ is a weak solution to problem (1.1). Assume that 
$v\in W_2^{2,1}(\Omega^T)$ and $u\in C^{1/2,1/4}(\Omega^T)$. Assume that 
${[v_\varphi^{(1)}(0)]^2\over r}$, ${\chi^{(1)}(0)\over r}\in L_2(\Omega)$.
Then
$$\eqal{
&\|v^{'(1)}\|_{V_2^1(\Omega^T)}\le\varphi(d_1,d_2,1/r_0)\bigg[1+
\bigg\|{[v_\varphi^{(1)}(0)]^2\over r}\bigg\|_{L_2(\Omega)}\cr
&\quad+\bigg\|{\chi^{(1)}(0)\over r}\bigg\|_{L_2(\Omega)}\bigg].\cr}
\leqno(5.1)
$$

\noindent
In view of (2.24) $\chi^{(2)}$ is a solution to the problem
$$\eqal{
&\chi_{,t}^{(2)}+v\cdot\nabla\chi^{(2)}-{v_r\over r}\chi^{(2)}-v\cdot\nabla
\zeta_2\chi\cr
&\quad-\nu\bigg[\bigg(r\bigg({\chi^{(2)}\over r}\bigg)_{,r}\bigg)_{,r}+
\chi_{,zz}^{(2)}+2\bigg({\chi^{(2)}\over r}\bigg)_{,r}\bigg]\cr
&\quad+\nu\bigg[{\chi\over r}\zeta_{2,r}^2+2\bigg({\chi\over r}\bigg)_{,r}
\zeta_{2,r}^2-\chi\zeta_{2,rr}^2\bigg]=
{2v_\varphi^{(2)}v_{\varphi,z}^{(2)}\over r},\cr
&\chi^{(2)}|_{S_1}=0,\ \ \chi^{(2)}|_{t=0}=\chi^{(2)}(0),\cr
&\chi^{(2)}|_{S_2}\ \ {\rm satisfic\ periodic\ boundary\ conditions\ with\ 
respect\ to}\ z.\cr}
\leqno(5.2)
$$

\proclaim Lemma 5.2. 
Assume that $v$ is a weak solution to problem (1.1). Assume that 
$v_\varphi(0)\in L_{49/18}(\Omega_{\zeta_2})$, $\chi^{(2)}(0)\in L_2(\Omega)$.
Then
$$
\|\chi\|_{V_2^0(\Omega_{\bar\zeta_2}^T)}\le c(1/r_0,d_1,d_2)(d_1^2+
\|\chi^{(2)}(0)\|_{L_2(\Omega)}),
\leqno(5.3)
$$
where $\Omega_{\bar\zeta_2}=\{x\in\Omega:\ \zeta_2(r)=1\}$.

\Proof 
Multiplying (5.2) by ${\chi^{(2)}\over r^2}$, integrating over $\Omega$ and 
using the boundary conditions yields
$$\eqal{
&{1\over2}{d\over dt}\intop_\Omega\bigg|{\chi^{(2)}\over r}\bigg|^2dx+
\nu\intop_\Omega\bigg|\nabla{\chi^{(2)}\over r}\bigg|^2dx
=\intop_\Omega v\cdot\nabla\zeta_2\chi{\chi^{(2)}\over r^2}dx\cr
&\quad-\nu\intop_\Omega
\bigg[{\chi\over r}\zeta_{2,r}^2+2\bigg({\chi\over r}\bigg)_{,r}\zeta_{2,r}^2-
\chi\zeta_{2,rr}^2\bigg]{\chi^{(2)}\over r^2}dx\cr
&\quad+2\intop_\Omega{v_\varphi^{(2)}v_{\varphi,z}^{(2)}\over r}
{\chi^{(2)}\over r^2}dx\cr}
\leqno(5.4)
$$
The first term on the r.h.s. is estimated by
$$
\varepsilon\bigg\|{\chi^{(2)}\over r}\bigg\|_{L_6(\Omega)}^2+
c(1/\varepsilon,1/r_0)\|v_r\|_{L_2(\Omega)}^2
\bigg\|{\chi\over r}\bigg\|_{L_3(\Omega_{\zeta_{2,r}})}^2,
$$
where $\Omega_{\zeta_{2,r}}=\Omega\cap\supp\zeta_{2,r}$
\goodbreak

\noindent
The second term on the r.h.s. of (5.4) we express in the form
$$
-\nu\intop_\Omega\bigg({\chi^2\over r^3}\zeta_{2,r}^2-{\chi^2\over r^2}
\zeta_{2,rr}^2\bigg)\zeta_2^2dx-2\nu\intop_\Omega
\bigg({\chi\over r}\bigg)_{,r}{\chi\over r}{1\over r}\zeta_{2,r}^2\zeta_2^2dx
\equiv I.
$$
Integrating by parts the second integral in $I$ takes the form
$$
\nu\intop_\Omega{\chi^2\over r^2}(\zeta_2^2\zeta_{2,r}^2)_{,r}drdz.
$$
Hence,
$$
|I|\le c(1/r_0)\intop_\Omega\chi^2dx.
$$
Finally, the last integral on the r.h.s. of (5.4) is estimated by
$$
\varepsilon\intop_\Omega\bigg|\nabla{\chi^{(2)}\over r}\bigg|^2dx+
c(1/\varepsilon)\intop_\Omega{|v_\varphi^{(2)}|^4\over r^4}dx.
$$
Using the above estimates in (5.4), assuming that $\varepsilon$ is sufficiently 
small and integrating the result with respect to time we obtain
$$\eqal{
&\bigg\|{\chi^{(2)}\over r}\bigg\|_{V_2^0(\Omega^t)}^2\le c(1/r_0)d_1^2
\bigg\|{\chi\over r}\bigg\|_{L_2(0,t;L_3(\Omega_{\zeta_{2,r}}))}^2\cr
&\quad+c(1/r_0)d_1^2+c(1/r_0)\intop_{\Omega^t}|v_\varphi^{(2)}|^4dxdt+
\bigg\|{\chi^{(2)}(0)\over r}\bigg\|_{L_2(\Omega)}^2,\cr}
\leqno(5.5)
$$
where estimate (2.21) is employed.
\goodbreak

\noindent
In view of (2.26) and the interpolation
$$\eqal{
&\bigg\|{\chi\over r}\bigg\|_{L_2(0,T;L_3(\Omega_{\zeta_{2,r}}))}^2\le
\varepsilon\bigg\|\nabla{\chi\over r}\bigg\|_{L_2(\Omega_{\zeta_2}^t)}^2+
c(1/\varepsilon)\bigg\|{\chi\over r}\bigg\|_{L_2(\Omega_{\zeta_2}^t)}^2\cr
&\le\varepsilon\bigg\|\nabla{\chi\over r}\bigg\|_{L_2(\Omega_{\zeta_2}^t)}^2+
c(1/\varepsilon,1/r_0)d_1^2,(1+d_2^2),\cr}
$$
we obtain from (5.5) the inequality
$$
\bigg\|{\chi\over r}\bigg\|_{V_2^0(\Omega_{\bar\zeta_2}^t)}^2\le
\varepsilon\bigg\|{\chi\over r}\bigg\|_{V_2^0(\Omega_{\zeta_2}^t)}^2+
c(1/\varepsilon,1/r_0,d_1,d_2)d_1^2
+\bigg\|{\chi^{(2)}(0)\over r}\bigg\|_{L_2(\Omega)}^2.
\leqno(5.6)
$$
By the local iteration argument (see [LSU, Ch. 4, Sect. 10]) inequality 
(5.6) for $\varepsilon<1$ implies (5.3). This concludes the proof.

In view of (5.3) we have the estimate for 
$\chi\in V_2^0(\Omega_{\bar\zeta_2}^T)$, where 
$\Omega_{\bar\zeta_2}=\{x\in\Omega:\ \zeta_2(r)=1\}$.
Take a function $\zeta_3=\zeta_3(r)$ such that 
$\supp\zeta_3\subset\Omega_{\bar\zeta_2}$.

\noindent
Now we consider the problem
$$\eqal{
&v_{r,z}-v_{z,r}=\chi,\cr
&v_{r,r}+v_{z,z}+{v_r\over r}=0,\cr}
\leqno(5.7)
$$
in $\Omega_{\bar\zeta_2}$. Hence, to formulate the problem more precisely 
we multiply (5.7) by $\zeta_3$ and introduce the notation
$$
v^{(3)}=v\zeta_3,\quad \chi^{(3)}=\chi\zeta_3.
$$
Then (5.7) takes the form
$$\eqal{
&v_{r,z}^{(3)}-v_{z,r}^{(3)}=\chi^{(3)}-v_z\zeta_{3,r}\quad &{\rm in}\ \ 
\Omega_{\bar\zeta_2},\cr
&v_{r,r}^{(3)}+v_{z,z}^{(3)}=-{v^{(3)}\over r}+v_r\zeta_{3,r}\quad &{\rm in}
\ \ \Omega_{\bar\zeta_2},\cr
&v^{(3)}\cdot\bar n|_{S_1}=0,\ \ v^{(3)}|_{S_2}\ \ {\rm satisfies\ periodic\ 
conditions,}\cr
&v^{(3)}|_{r=r_0}=0.\cr}
\leqno(5.8)
$$
Since the r.h.s. of (5.8) belongs to $V_2^0(\Omega_{\bar\zeta_2}^T)$ 
we obtain the estimate
$$
\|v^{'(3)}\|_{V_2^1(\Omega_{\bar\zeta_2}^T)}\le c
(\|\chi^{(3)}\|_{V_2^0(\Omega_{\bar\zeta_2}^T)}+d_1).
\leqno(5.9)
$$
Hence, from (5.1) and (5.9) we obtain the estimate
$$\eqal{
&\|v'\|_{V_2^1(\Omega^T)}\le\varphi(d_1,d_2,1/r_0)\bigg[d_1^2+
\bigg\|{|v_\varphi^{(1)}(0)|^2\over r}\bigg\|_{L_2(\Omega)}\cr
&\quad+\bigg\|{\chi^{(1)}(0)\over r}\bigg\|_{L_2(\Omega)}+
\|\chi^{(2)}(0)\|_{L_2(\Omega)}\bigg]
\equiv d_3,\cr}
\leqno(5.10)
$$
where $v'=(v_r,v_z)$.

\noindent
From [Z2, proof of Lemma 3.7] and (5.10) we have
$$
\|v'\|_{L_{10}(\Omega^T)}+\|\nabla v'\|_{L_{10/3}(\Omega^T)}\le cd_3\equiv A_0.
\leqno(5.11)
$$

\proclaim Theorem 5.3. (Global existence) (Proof of Theorem A) 
Let the Assumptions from Section 1 hold. Then there exists a solution to 
problem (1.1) such that $v\in W_2^{2,1}(\Omega^T)$, 
$\nabla p\in L_2(\Omega^T)$, $T\in R_+$ and
$$
\|v\|_{W_2^{2,1}(\Omega^t)}+\|\nabla p\|_{L_2(\Omega^t)}\le A_1+A_2+A_3
\equiv\alpha,
\leqno(5.12)
$$
where $t\in[(k-1)T_0,kT_0]$, $k\in\N$, $T_0$ is such that
$$
c_*T_0^{1/2}\alpha\le1
\leqno(5.13)
$$
and $c_*$ is introduced in (3.3).

\Proof 
Assume that $v(0)\in H^1(\Omega)$ is given. Then Theorem 3.1 implies existence 
of local solution $v\in W_2^{2,1}(\Omega^{T_*})$, where $T_*$ satisfies (3.3). 
By imbedding $v\in L_{10}(\Omega^{T_*})$ so [Z4] yields
$$
\|u\|_{C(0,T_*;C^{1/2}(\Omega))}\le d_4,
\leqno(5.14)
$$
where $d_4$ does not depend on the local solution. Then (5.10) implies that
$$
\|v'(T_*)\|_{H^1(\Omega)}\le d_3\le A_1+A_2\equiv A_0,
\leqno(5.15)
$$
where $A_0$ does not depend on $T_*$.

\noindent
Hence (5.15) and (4.1) imply
$$
\|v(T_*)\|_{H^1(\Omega)}\le d_3+d_5\le A_1+A_2+A_3\equiv\alpha.
\leqno(5.16)
$$
Choosing $T_0<T_*$ satisfying (5.13) we can continue the above considerations 
in the interval $[T_0,2T_0]$. In this way we show that
$$
\|v(kT_0)\|_{H^1(\Omega)}\le\alpha.
\leqno(5.17)
$$
This concludes the proof.

\section{References}

\item{BIN.} Besov, O. V.; Il'in, V. P.; Nikolskij, S. M.: Integral 
representation of functions and theorems of imbedding, Nauka, Moscow 1975.

\item{L.} Ladyzhenskaya, O. A.: The mathematical theory of viscous 
incompressible flow, Nauka, Moscow (in Russian).

\item{LSU.} Ladyzhenskaya, O. A.; Solonnikov, V. A.; Uraltseva, N. N.: 
Linear and quasilinear equations of parabolic type, Nauka, Moscow 1967 
(in Russian).

\item{Z1.} Zaj\c aczkowski, W. M.: A priori estimate for axially symmetric 
solutions to the Navier-Stokes equations near the axis of symmetry

\item{Z2.} Zaj\c aczkowski, W. M.: Global special regular solutions to the 
Navier\--Stokes equations in a cylindrical domain without the axis of 
symmetry, Top. Meth. Nonlin. Anal. 24 (2004), 69--105.

\item{Z3.} Zaj\c aczkowski, W. M.: Global special regular solutions to the 
Navier\--Stokes equations in a cylindrical domain under boundary slip 
condition, Gakuto Intern. Ser., Math. Sc. Appl. 21 (2004), pp. 188.

\item{Z4.} Zaj\c aczkowski, W. M.: The H\"older continuity of the swirl for 
the Navier-Stokes motions

\item{Z5.} Zaj\c aczkowski, W. M.: On global regular solutions to the 
Navier-Stokes equations in cylindrical domains, Top. Meth. Nonlin. Anal. 
37 (2011), 55--85.

\item{S.} Solonnikov, V. A.: Estimates of the solutions of a nonstationary 
linearized system of the Navier-Stokes equations, Trudy Mat. Inst. Steklov 
70 (1964), 213--317; English transl. Amer. Math. Soc. Trans. Ser. 2, 65 
(1967), 51--137.

\item{ZZ.} Zadrzy\'nska, E.; Zaj\c aczkowski, W. M.: Nonstationary Stokes 
system in Sobolev spaces.

\item{NP1.} Neustupa, J.; Pokorny, M.: Axisymmetric flow of Navier-Stokes fluid 
in the whole space with non-zero angular velocity component, Math. Bohem. 
126 (2001), 469--481.

\item{NP2.} Neustupa, J.; Pokorny, M.: An interior regularity criterion 
for an\break 
axially symmetric suitable weak solution to the Navier-Stokes equation, 
J. Math. Fluid Mech. 2 (2000), 381--396.

\item{P.} Pokorny, M.: A regularity criterion for the angular velocity 
component in the case of axisymmetric Navier-Stokes equations, in Elliptic 
and Parabolic Problems (Rolduc/Gaeta, 2001), World Scientific, River Edge, 
NJ (2002), 233--242.

\bye